\documentclass[a4paper, 10pt, twoside]{article}
\usepackage[english]{babel}
\usepackage{amsmath, amssymb,theorem}
\usepackage[all]{xy}
\usepackage{epsfig}
\UseComputerModernTips
\title{Operads up to Homotopy and\\ Deformations of Operad Maps
}
\author{Pepijn van der Laan}
\makeatletter
\renewcommand\subsection{\@startsection{subsection}{2}{\z@}%
                                     {-2.3ex\@plus -1ex \@minus -.2ex}%
                                     {1ex \@plus .2ex}%
                                     {\normalfont\bfseries\centering}}
\makeatother

\newcommand{\note}{\textbf}

\newcommand{\NN}{\mathbb{N}}
\newcommand{\ZZ}{\mathbb{Z}}

\newcommand{\Gh}{\mathfrak{h}}

\newcommand{\Gg}{\mathfrak{g}}

\renewcommand{\phi}{\varphi}

\newcommand{\eps}{\varepsilon}

\newcommand{\pd}{\partial}

\newcommand{\colim}[1]{\underset{#1}{\text{colim}}}
\newcommand{\End}{\text{End}}
\newcommand{\Aut}{\text{Aut}}
\newcommand{\id}{\text{id}}

\newcommand{\Hom}{\text{Hom}}
\renewcommand{\vert}{\mathbf{v}}
\newcommand{\edge}{\mathbf{e}}
\newcommand{\leg}{\mathbf{l}}
\newcommand{\Fin}{\textsf{Fin}}

\newcommand{\gVect}{\textsf{-gVect}}
\newcommand{\dgVect}{\textsf{-dgVect}}

\newcommand{\Set}{\textsf{Set}}
\newcommand{\Coopd}{\textsf{Coopd}}
\newcommand{\Opd}{\textsf{Opd}}

\newcommand{\Sh}{\text{Sh}}
\newcommand{\fr}{\frac}

\newcommand{\surjection}{\twoheadrightarrow}
\newcommand{\be}{\begin{equation}}
\newcommand{\ee}{\end{equation}}
\theoremheaderfont{\scshape}
\theorembodyfont{\itshape}
\theoremstyle{change}
\newtheorem{Tm}{Theorem}[section]
\newtheorem{Pp}[Tm]{Proposition}
\newtheorem{Lm}[Tm]{Lemma}
\newtheorem{Cr}[Tm]{Corollary}

\theorembodyfont{\rmfamily}
\newtheorem{Cv}[Tm]{Convention}
\newtheorem{Df}[Tm]{Definition}
\newtheorem{Ex}[Tm]{Example}
\newtheorem{Rm}[Tm]{Remark}

\newenvironment{Pf}{{\scshape Proof}}{\hspace*{\fill}{\scshape QED}\\
\par}
\begin{document}
\maketitle
%\chapter{Models and Deformation}

%\parskip=0cm
%\tableofcontents
%\parskip=.1cm

\section{Introduction}

\subsection{Abstract}
From the `cofree' cooperad $T'(A[-1])$ on a collection $A$ together with a 
differential, we construct an $L_\infty$-algebra structure on the total
space 
$\bigoplus_nA(n)$ that descends to coinvariants. We use this construction 
to define  an $L_\infty$-algebra controlling deformations of the
operad $P$ under $Q$ from a cofibrant resolution for an operad $Q$,
and an operad map $Q\longrightarrow P$. Starting from a diffent
cofibrant resolution one obtains a quasi isomorohic
$L_\infty$-algebra. 
This approach unifies Markl's cotangent cohomology of operads and the
approaches to deformation of $Q$-algebras by Balavoine, and
Kontsevich and Soibelman.
 
\subsection{Motivation}

We recall that the total space $\bigoplus_nP(n)$ of any
(in $k\dgVect$) operad $P$ has a natural Lie algebra structure that
descends to the coinvariants  $\bigoplus_nP(n)_{S_n}$. Since
bar/cobar duality identifies operad structures on the collection $P$
with differentials on the `cofree' cooperad $T'(P[-1])$, it is natural
to ask what 
happens on $P$ if we start from an arbitrary differential on $T'(P[-1])$. We call
the  structure obtained on $P$ an operad up to homotopy. The explanation of this
terminology in terms of Quillen model categories is postponed to
\cite{Pep:colopd}.
Our notion is manifestly different from the lax operads described by 
Brinkmeier \cite{Brin:Thesis}. In our context the compatibility with
the symmetric group actions is strict. This makes the notion more applicable.

We show that an operad up to homotopy $P$ gives rise to a natural 
$L_\infty$-algebra structure on the total space $\bigoplus_nP(n)$, that 
descends to coinvariants. This is the natural analogue of the Lie algebra
structure on the total space of on operad. Morphisms of operads up to
homotopy 
(i.e. morphisms of the corresponding free cooperads with differential) yield
morphisms between the corresponding $L_\infty$-algebras.

Applying the Lie algebra construction to the convolution operad $P^A$
of a cooperad $A$ and an operad $P$, solutions to the Maurer-Cartan
equation correspond to operad maps $B^*A\longrightarrow P$ from the
cobar construction $B^*A$ of $A$ to $P$. Now suppose that $B^*A$ is a
cofibrant resolution for an operad $Q$, then the Lie algebra
$\bigoplus_nP^A(n)_{S_n}$ controlls deformations of operads 
$Q\longrightarrow P$ under $Q$.  We show the Lie algebras for 
different cofibrant models are quasi isomorphic. We use the model 
category of operads and the functoriality of the $L_\infty$-structure 
on the total space to do this. 
We observe that this construction extends to more general models for
$Q$, as to include Markl's cotangent cohomology of operads, and the
approaches to deformations of $Q$-algebras by Balavoine, and
Kontsevich and Soibelman.

The material in this paper is intended to constitute the core of a
chapter in the author's dissertation. The author aims to include more
applications and examples there. 
Please note the paper is written in a way that allows generalisation to
coloured operads.

\subsection{Plan of the Paper}

Section 2 gives a short overview of the necessary standard results on
operads. Section 3 introduces operads up to homotopy and shows how to
construct $L_\infty$ algebras from these. We show how morphisms of
operads up to homotopy induce morphisms of
$L_\infty$-algebras. Section 5 defines a cohomology
theory using convolution operads. The $L_\infty$-algebra constructed
in Section 3 serves to show that the cohomology is independent of the
necessary choice of cofibrant resolution. Secion 5 continues to show
the relation the known approaches to deformations of algebras over an operad.
To fulfill these tasks we need the results we recall in Section 4.

\subsection{Acknowledgements}
This paper emerged in answer to a question posed to the author by Jeff
Smith: "Is there a deformation theory of operad maps?"
The author is grateful to both Ieke Moerdijk and Martin Markl for their
interest and for many illuminating discussions, to Jeff Smith for
the question, and to Alessandra Frabetti for her comments on a
previous draft.
\section{Operads and Cooperads}
This section is purely expositorial. We describe relevant structures
on operads and cooperads such as the Bar/Cobar construction and the
convolution operad.
\subsection{Definitions}
Throughout this paper, let $k$ be a field of characteristic 0. We work
in the category 
$k\dgVect$ of dg vector spaces over $k$ (by convention $\ZZ$-graded with a
differential of degree +1). This is a symmetric monoidal category with the
usual tensor product $\otimes$ and the symmetry $A\otimes B \ni a\otimes b
\mapsto (-1)^{|a||b|}b\otimes a\in A\otimes B$ for homogeneous elements $a$
and $b$ of degrees $|a|$ and $|b|$ is dg vector spaces $A$ and $B$.
We will make extensive use of the shift functor on $k\dgVect$, given by 
$V[p]^n = V^{n-p}$ with differential $d[p]^n = (-1)^p d^{n-p}$ (superscripts
denoting degrees here, not powers).
We have the usual Koszul sign convention. In particular, for graded maps $f$
and $g$ we write $(f\otimes g)(a\otimes b) = (-1)^{|g||a|}f(a)\otimes g(b)$.
These conventions allow us to hide most of the signs.

Let $\Fin_*$ be the category of pointed finite sets and bijections
preserving the basepoint. 
For $(X,x_0), (Y,y_0) \in \Fin_*$ and $x\in X-\{x_0\}$ define $Y\cup_xX$
to be the pushout (in $\Set$, the category of sets)
\[
\xymatrix{
\ast \ar[d]_{y_0} \ar[r]^{x} \ar@{}[dr]|>{\large \ulcorner}& X\ar[d] \\ 
Y\ar[r] & Y\cup_xX,
}
\]
with the basepoint $x_0$. 
%That is, as an element of $\Set$, it is the
%direct union of $X$ and $Y$ modulo the equivalence relation generated
%by $y_0\sim x$. 
Let $\tau\in\Aut(Y,y_0)$ and
$\sigma\in\Aut(X,x_0)$. The map $\sigma\circ_x \tau:X\cup_xY
\longrightarrow X\cup_{\tau(x)}Y$ in $\Fin_*$ is the unique dotted arrow
that makes the diagram below commute
\begin{equation}\label{eq:pushauto}
\xymatrix{ Y \ar[r]\ar[d]^{\sigma}& Y\cup_xX \ar@{.>}[d] 
& X \ar[l]\ar[d]^\tau \\
Y \ar[r] & Y\cup_{\tau(x)}X & \ar[l]X.
}
\end{equation}
A \note{collection} in is
a contravariant functor $C:\Fin_*\longrightarrow k\dgVect$. 
A \note{pseudo operad} is a collection together with a map
\[
\circ_x:C(X,x_0)\otimes C(Y,y_0) \longrightarrow C(X\cup_{x}Y,x_0),
\]
for each pair of
pointed sets $(X,x_0)$ and $(Y,y_0)$ and each $x\in X$ such that
$x\neq x_0$, such that
\begin{enumerate}
\item
With respect to the automorphisms of $X$ and $Y$ the map $\circ_x$
behaves as
\[
p\tau\circ_xq\sigma = (p\circ_{\tau^{-1}x} q)
(\tau\circ_{\tau^{-1}x}\sigma)
\]
(cf. equation (\ref{eq:pushauto})) for $p\in C(X), q\in C(Y)$ and $x\in
X-\{x_0\},\tau\in\Aut(X,x_0),\sigma\in\Aut(Y,y_0)$.
\item The squares
\begin{align*}
\xymatrix{ C(X)\otimes C(Y)\otimes C(Z)
\ar[r]^{\circ_x\otimes\id}\ar[d]^{(\id\otimes\circ_{x'})\circ
  (\id\otimes s)} & C(X\cup_xY)\otimes
C(Z) \ar[d]^{\circ_{x'}} \\
C(X\cup_{x'}Z)\otimes C(Y) \ar[r]^{\circ_x} & C(X\cup_x Y \cup_{x'} Z)}
\\
\xymatrix{ C(X)\otimes C(Y)\otimes C(Z)
\ar[r]^{\id\otimes\circ_y}\ar[d]_{\id\otimes \circ_{x}} & C(X)\otimes
C(Y\cup_yZ)\ar[d]^{\circ_{x}} \\
C(X\cup_{x}Y) \otimes C(Z) \ar[r]^{\circ_y} & C(X\cup_x Y \cup_y Z)
}
\end{align*}
%\xymatrix{
% C(Y)\otimes C(X) \ar[r]^{\sigma\otimes\tau}\ar[d]^{\circ_x} &
% C(Y)\otimes C(X) \ar[d]^{\circ_{\sigma(x)}}\\
% C(Y\cup_xX) \ar[r]^{\tau\circ_x\sigma} &C(Y\cup_{\sigma(x)} X)
%}
commute for $y\in Y$ and $x,x'\in X$. We omit the basepoints from the
notation and the map $s$ is the symmetry of the tensor product.
\end{enumerate}
An \note{operad}
is a pseudo operad together with a map $\id:I\longrightarrow
C(\{x,x_0\})$ for a two element set $\{x,x_0\}$, such that $\id$ is a
right/left identity with respect to any well defined right/left composition. 
Together with the natural notion of morphism this defines the
categories of collections, pseudo operads and operads.  

In some cases (though not always) a different description is convenient. 
The description is obtained by the inclusion of
the symmetric groupoid as a skeleton in $\Fin_*$ with objects
$\mathbf{n} = \{0,\ldots,n\}$ and 0 as basepoint. We write $P(n) =
P(\mathbf{n})$. At times we will use
the equivalent description of operads in terms of the maps 
\begin{equation}\label{eq:gamma}
\gamma: P(n)\otimes_{S_n}(P(k_1)\otimes\ldots\otimes P(k_n)
\longrightarrow P(k_1+\ldots + k_n)
\end{equation}
and an identity map in $P(1)$. Here $\gamma(p,q_1,\ldots,q_n) =
(\cdots((p\circ_nq_n)\circ_{n-1}q_{n-1})\ldots \circ_1q_1)$ for $p\in
P(n)$, and $q_i\in P(k_i)$.

A \note{(pseudo) cooperad} is a (pseudo) operad in the opposite category 
$(k\dgVect)^{\text{op}}$. That is, a covariant(!) functor
$\Fin_*\longrightarrow k\dgVect$ satisfying the dual axioms. 
Here the directions of arrows is inverted, the coinvariants change to
invariants, sums change into products, and the right action 
of $S_n$ changes to a left action. For convenience of the reader 
we twist this action by an inverse and obtain a right
action. Therefore we do not have to distinguish between left and right
collections.

\subsection{Graphs and Trees}

This paragraph recalls the approach of
Getzler and Kapranov \cite{GtzKap:Mod} on graphs. 
A \note{graph} is a finite set of \note{flags}
(or 
\note{half edges}) together with both an equivalence relation $\sim$,
and an automorphism $\sigma$  of order $\leq 2$. The equivalence
classes are called \note{vertices}, the orbits of length 2
\note{(internal) edges}, and the orbits of length 1 \note{external
edges} or \note{legs}. We denote the vertices of a graph $g$ by
$\vert(g)$, its (internal) edges by $\edge(g)$ and its legs by
$\leg(g)$, especially we denote the half edges in a vertex $v$ by
$\leg(v)$. The flags in a vertex are also called its
\note{legs}. The vertex containing a flag $f$ is called the vertex of $f$.
To picture a graph, draw a node for each vertex with outgoing edges
labeled by its flags. For each orbit of
$\sigma$ of length 2 (edge) draw a line connecting the two flags. For
each orbit of length one this leaves an external edge. 
A \note{tree} is a connected graph $t$ such that $|\vert(t)| -
|\edge(t)| =1$. A \note{rooted tree} is a tree together with a base
point in the 
set of external edges. Note that on a rooted tree, the base point 
induces a canonical base point on the legs of any vertex. We denote by
$T(X)$ the set of rooted trees with external edges labeled by the
pointed set $X$. (An element $t$ of $T(X)$ is a tree $t$ together with
a basepoint preserving bijection $X\longrightarrow \leg(t)$.)
A \note{morphism} of graphs is a morphism $\phi$ of the set of flags,
such that for all flags $f,f'$ both  $\phi(\sigma f) = \sigma \phi(f)$, and
$\phi(f)\sim\phi(f')$ iff $f\sim f'$ (where $\sigma$ and $\sim$ are
the defining automorphisms and equivalence relations). If the set of
external edges has a basepoint, we assume that morphisms preserve
the base point. The \note{group of automorphisms} of a graph $g$ is denoted
$\Aut(g)$. Note that an automorphism defines an isomorphism of
flags, an automorphism of edges, an automorphism of external edges,
and an isomorphism of vertices. Moreover it defines an isomorphism of
the legs of a vertex to the legs of its image under the map on vertices. 
For two rooted trees $s$ and $t$ and $x\in \leg(t)$ define $s\circ_xt$
as the rooted tree obtained from $s$ and $t$ by grafting the root of
$s$ on leg $x$ of $t$.
 
\subsection{Free Operads and `Cofree' Cooperads}

 The forgetful functor from pseudo operads to collections in
$k\dgVect$ has a left adjoint 
$T$, which is the \note{free pseudo operad} functor. We can give
an explicit formula for ${T}$. Let
\begin{equation}\label{eq:Ctdef}
C(t) = \left(\bigotimes_{v\in \vert(t)}C(\leg(v))\right)_{\Aut(t)},
\end{equation}
for any rooted tree $t$. The action of $\Aut(t)$ is by permutation of the
tensor factors according to the permutation of vertices, by the associated
isomorphism of $\leg(v)\longrightarrow \leg(\sigma(v))$ action inside the tensor
factor associated to $v$, and by the composition of the labeling of
$\leg(t)$ by the induced isomorphism of $\leg(t)$. 
Define the underlying collection ${T}C$ of the free pseudo operad by
\[
{T}C(X) = \colim{T(X)} (C(t)) \equiv \bigoplus_{t\in \bar{T}(X)} C(t),
\]
where the colimit is over the groupoid with objects $t\in T(X)$ and
isomorphisms of trees as maps, and the sum is over the isomorphism classes
of rooted trees $\bar T(X)$ in $T(X)$. 
The pseudo operad structure on $TC$ is given by grafting trees. 
%For trees $s$ and $t$, and an external edge $x$ of $s$, the composition 
%$\circ_x:TC(t)\otimes TC(s)\longrightarrow TC(t\circ_x s)$ is the obvious
%the map that rearranges the tensor factors and takes the quotient with
%respect to $\Aut(t\circ_x s)$. Associativity and equivariance
%are immediate. 
The forgetful functor from operads to pseudo operads has an left
adjoint. This map is given by adjoining a unit with respect to
composition in $P(X)$ for $X$ such that $|X|=2$. By definition the
\note{free operad functor} $T^+$ is the composition of both left adjoints.

We now define a cooperad that is cofree with respect to a
restricted class of cooperads. Therefore we use the term `cofree' when
referring to this cooperad. The `cofree' pseudo cooperad functor
${T}'$ is defined by 
\[
{T}'C(X,x_0) = \bigoplus_{t\in \bar T(X)} \left(
\bigotimes_{v\in\vert(t)}C(\leg(v)) \right)^{\Aut(t)}.
\]
The pseudo cooperad structure is determined by cutting edges.
%We have an operation $\circ_e^*$, the dual to the operation $\circ_x$
%grafting $s$ and $t$, where $e$ is the newly created edge in
%$s\circ_x t$: Just seperate the tensor factors according to
%$v\in\vert(t)$ or $v\in \vert(s)$. 
%This gives a well defined map
%\[
%\left(\bigotimes_{v\in\vert(s\circ_xt)}C(\leg(v))
%\right)^{\Aut(s\circ_xt)}\longrightarrow 
%\left(\bigotimes_{v\in\vert(s)}C(\leg(v))
%\right)^{\Aut(s)}\otimes
%\left(\bigotimes_{v\in\vert(t)}C(\leg(v)) \right)^{\Aut(t)}
%\]
%that satisfies the desired property. 
Composition of the functor ${T}'$ with the right adjoint to the
forgetful functor from cooperads to pseudo cooperads defines the
\note{`cofree' cooperad functor} $(T')^+$. 

Recall the  monoidal structure $(\square,I,a,i_r,i_l)$ on the category
of collections in  $k\dgVect$ from Smirnov \cite{Smir:box} (Shnider
and Van Osdol \cite{ShnVOs:Opd} is a more accessible
reference):
\[
C\square D (X)= \bigoplus_Y\bigoplus_{f:X\surjection Y}
C(Y)\otimes_{\Aut(Y)}\left( \bigotimes_{y\in Y-y_0}
D(f^{-1}(y)\cup\{y\})\right)  
\]
where the sum is over one representing pointed sets $Y$  for each 
isomorphism class, and where $y$ is the basepoint of
$f^{-1}(y)\cup\{y\}$. Operads are algebras in this monoidal
category. 
Getzler and Jones \cite{GetzJon:Opd}
define a cooperad as a coalgebra with respect to this monoidal structure.
Since char$(k) = 0$, the natural isomorphism from coinvariants to
invariants identifies both notions of cooperad.

%
% Convolution
%

\subsection{Convolution Operad}

For convenience the lemma below uses the description of operads and
cooperads with an explicit choice of coordinates. That is, we use
$\gamma$ (cf. Formula \ref{eq:gamma}), and for cooperads the dual map
$\gamma^*$.  
Let $C$ and $P$ be collections.  
Consider the collection $P^C(X) =
\Hom(C(X),P(X))$, with the natural action 
(i.e. $(\phi\sigma)(c) := (\phi(c\sigma^{-1}))\sigma$ for $\sigma\in\Aut(X)$)
and the natural differential (i.e. $d\phi := \phi\circ d -d\circ\phi$).
 
\begin{Lm}[Berger and Moerdijk \cite{BerMoer:Model}] \label{Lm:convol}
Let $P$ be an operad and let $C$ be a cooperad. The collection $P^C$
enjoys an operad structure given by 
\[
\gamma(\phi; \psi_1,\ldots, \psi_n) 
= \gamma_P \circ (\phi\otimes (\psi\otimes\ldots\otimes \psi_n)) \circ
\gamma^*_C, 
\]
where we restrict $\gamma^*$ to the suitable summand for homogeneous
(wrt arity) $\phi$ and $\psi_i$. 
The unit is given by the composition $u\circ\eps$ of the
counit $\eps$ of $C$ with the identity $u$ of $P$.
\end{Lm}
%\begin{Pf}
%For  $\phi\in P^C(n)$ the differential is defined by $ \phi\circ d -
%d\circ\phi$,
%in terms of the differentials on $C$ and $P$. 
%The identity is defined as $u\circ \eps$, where $u$ is the unit of $P$ and 
%$\eps$ is the counit of $C$. 
%Associativity is directly clear from coassociativity of $\gamma^*_C$ and
%associativity of $\gamma_P$. It suffices to check the commutation
%relation of $\gamma$ with the action of $\Aut(\mathbf{n})$. Let
%$\phi\in P^C(\mathbf{n})$ and for $k\in\{1,\ldots,n\}$, let
%$\phi_k\in P^C(Y_k)$ for some pointed sets $Y_k$. We use a generalised
%Sweedler notation for the cocomposition in $C$ and compute 
%\[
%\begin{split}
%\sum_{(c)}\gamma_P(\phi\sigma(c_0);&
%\phi_{\sigma 1}(c_1),\ldots,\phi_{\sigma n}(c_n)) \\
%& = \sum_{(c)} \gamma_P(\phi(c_0\sigma^{-1})\sigma;
%\phi_{\sigma 1}(c_{1}),\ldots,\phi_{\sigma n}(c_{n})) \\
%& = \sum_{(c)} \gamma_P(\phi(c_0\sigma^{-1}); 
%\phi_{1}(c_{\sigma^{-1} 1}),\ldots,\phi_n(c_{\sigma^{-1} n})) \\
%& = \sum_{(c)}\gamma(\phi(c_0);\phi_1(c_1),\ldots,\phi_n(c_n)).
%\end{split}
%\]
%The first equality is the definition of the action, the second
%equality is by coinvariance of $\gamma_P$ the final equality is by
%invariance of $\gamma^*_C$.
%\end{Pf}
The operad defined in Lemma \ref{Lm:convol} is called the \note{convolution
operad} of $C$ and $P$. 

%
% Bar Cobar
%
\subsection{Bar and Cobar Construction}

An operad $P$ is \note{augmented} if there is a projection of operads
$\eps: P \longrightarrow k\cdot\id$
such that the unit $k\cdot\id\longrightarrow P(1)$ composed with $\eps$ is the
identity on $k\cdot\id$.  The kernel of $\eps$ is the augmentation ideal of
$P$. The categories of augmented operads and pseudo operads are isomorphic. 
The dual notion is a \note{coaugmented cooperad}.

\begin{Cv}\label{Cv:speudoopd}
In the sequel we use augmented operads exclusively.
Likewise for cooperads. We denote by $T$
(resp. $T'$) the (co)free (co)operad in the (co)augmented
category. 
\end{Cv}
Let $C$ be a graded collection. A \note{differential} on the `cofree'
cooperad $T'C$ is a coderivation $\pd$ of degree
$+1$  that satisfies $\pd^2=0$.
A differential on the free operad $TC$ is a derivation $\pd$ of degree
$+1$  that satisfies $\pd^2=0$.
Recall that the `cofree' cooperad is graded by the number of vertices
in the trees. 
\begin{Tm}[Getzler and Jones \cite{GetzJon:Opd}]\label{Tm:OpdCoder}
Let $P$ be a graded collection. There is a 1-1 correspondence between
operad structures on $P$ and differentials on $T'(P[-1])$ that are of
degree $\geq -1$ in vertices.

Dually, there is a 1-1 correspondence between cooperad structures on
$P$ and differentials on $T(P[1])$ that are of degree $\leq +1$ in vertices.
\end{Tm}
Let $P$ be an operad. Define the \note{bar construction} $B^*P$ to be the
corresponding cooperad $(T'(P[1]),\pd)$ of Theorem \ref{Tm:OpdCoder}.
Let $C$ be a cooperad. Define the \note{cobar construction} $BC$ to be the
corresponding operad $(T(P[1]),\pd)$ of Theorem \ref{Tm:OpdCoder}.
These constructions define two functors 
\[
B:\Opd\longrightarrow \Coopd \qquad \text{ and }\qquad B^*:\Coopd
\longrightarrow \Opd.
\] 
\begin{Tm}[Ginzburg and Kapranov \cite{GtzKap:Mod}]
Let $P$ be an operad. The natural projection $B^*BP \longrightarrow P$
is a quasi isomorphism of operads (i.e. the induced morphisms in
cohomology are isomorphisms).
Dually, let $C$ be a cooperad.
The natural inclusion $C \longrightarrow BB^*C$ is a quasi
isomorphism of cooperads.
\end{Tm}

%
% O P E R A D S   U P T O   H T P Y
%
\section{Operads up to homotopy}
%
% Definitions
%
This section relaxes the operad
axioms and obtains thus the definition of an operad up to
homotopy. It shows an operad up to homotopy gives rise to an
$L_\infty$-algebra, and that a morphism of operads up to homotopy
gives rise to an $L_\infty$-morphism. 
The author will discuss the characterization of
operads up to homotopy in terms of homotopy algebras for the Koszul
self-dual $\NN$-colored Koszul operad of non-$\Sigma$ operads in
\cite{Pep:colopd}.

\subsection{Definitions}

Following Convention \ref{Cv:speudoopd} we use (co)augmented
(co)operads throughout the sequel of this paper.
\begin{Df}
An \note{operad up to homotopy} $P$ is a collection $P$ in $k\gVect$
together with a differential $\pd$ on the tree cooperad $T'(P[1])$.
Dually, a \note{cooperad up to homotopy} $C$ is a collection $C$ in $k\gVect$
together with a differential $\pd$ on the free operad $T(C[-1])$. 
Analoguous to the $A_\infty$-case one could discuss several notions of
units for operads up to homotopy. We ignore this point since it is
beyond our aims.
\end{Df}

\begin{Rm}
In the sequel of this section we restrict to operads up to homotopy.
The dual statements are readily deduced. 
\end{Rm}

\begin{Pp}[Markl \cite{Mar:Model}]
Let $C$ be a collection.
There is a 1-1 correspondence between coderivations of $T'C$ and
collection morphisms $T'C\longrightarrow C$.
\end{Pp}

The differential on $(T'P[-1],\pd)$ defines for each rooted tree $t$ an
operation
\[
\circ_t:P[-1](t)\longrightarrow P[-1](\leg(t))
\]
(cf. Equation \ref{eq:Ctdef}), that preserves both internal and
collection degree. Conversely, $\pd$ is completely determined by the
operations $\circ_t$.
Part of the structure is an internal differential $d$ on the collection
$P$ (corresponding to trees with one vertex).
The condition on $\pd^2 = 0$ on the differential is equivalent to a
sequence of relations on these 
operations. For each rooted tree $t$, we obtain a relation of the form
\begin{equation}\label{eq:pdsquare}
\sum_{s\subset t}\pm(\circ_{t/s})\circ(\circ_{s}) = 0,
\end{equation}
where the sum is over connected subtrees $s$ of $t$ and $t/s$ is the
tree obtained  from $t$ by contracting the subtree $s$ to a point.
The signs involved are induced by a choice of ordering on the vertices
of the trees $t$ and $s$ and the Koszul convention. 
If $P$ is an operad up to 
homotopy, then the cohomology $H^*P$ with respect to the internal differential
$d$ induced by $\pd$ is a graded operad.
\begin{Df}
An operad
up to homotopy is called \note{minimal} if the composition of  $\pd|_{P[1]}$
with the projection to the cogenerators of $T'(P[1])$ is zero. That
is, if the differential on the collection $P$ as induced by $\pd$
vanishes, or equivalently such that the differential on $T'(P[1])$ is of 
degree $\leq -1$ with respect to vertices.
An operad up to homotopy is called \note{strict} if the differential
$\pd$ is of degree $\geq -1$ in the number of vertices. 
Through Theorem \ref{Tm:OpdCoder} there is an obvious 1-1
correspondence between operads and strict operads up to homotopy. 
\end{Df}
%\begin{Rm}
%Minimality does not imply strictness. In the case of $L_\infty$-algebras
%we have
%natural examples of this phenomenon in the case of Courant algebroids over
%a point 
%that form a Lie algebra together with a non-trivial operation $l_3$ (cf.
%Roythenberg \cite{}). In the operad case, we encounter examples as minimal
%models for non-Koszul quadratic operads.
%\end{Rm}

\subsection{Planar Rooted Trees}
In order to show that for any operad up to homotopy $P$ the total space
$\pmb{\oplus}P$ is an $L_\infty$-algebra we need planar trees.
These will serve to keep track of the $\Aut(\leg(v))$-actions, for
$v\in\vert(t)$ for a rooted tree $t$.
\begin{Cv}
Let the functor
$T'$ from $k\dgVect$ to cooperads denote the `cofree' coassociative
coalgebra. On a vector space $V$ we have $T'V=\bigoplus_{n>0}V^{\otimes n}$.
Confusion with the `cofree' cooperad functor on
collections should not be possible.
Let $S'$ denote the `cofree' non-counital cocommutative 
coalgebra functor, on a vector space $V$ given by $S'V =
\bigoplus_{n>0}(V^{\otimes n})^{S_n}$. 
Introduce two functors $\pmb{\oplus}$ and $\pmb{\oplus}_S$ from
collections to $k\dgVect$. On an object $C$ these are given by
\[
\pmb{\oplus}(C) = \bigoplus_{n\in\NN} C(n)\qquad \text{and}\qquad
\left.\pmb{\oplus}\right._S(C) =\bigoplus_{n\in\NN} (C(n))_{S_n}.
\]
\end{Cv}
For $n\in\NN$, denote by $\mathbf{n}$ the pointed set $\{0,\ldots,n\}$
with basepoint $0$.
A \note{planar rooted tree} is a rooted tree together with  an
isomorphism of pointed sets $\mathbf{n}_v \longrightarrow \leg(v)$ for each
$v\in \vert(t)$, that sends $0$ to the
base point, where $n_v=|\leg(v)|-1$. 
%For each isomorphism class of rooted trees we fix
%a representing element. 
The planar structure induces a linear ordering
on the vertices (we write $\vert(t) = \{v_1,\ldots,v_{|\vert(t)|}\}$)   
and a labeling of $\leg(t)$ by the pointed set 
$\mathbf{n}_t = \{0,\ldots,|\leg(t)|-1\}$.
To see this, draw a planar tree with half edges of $v\in\vert(t)$ from
left to right with respect to the ordering of
$\leg(v)$.
%
%Given $m$, for each isomorphism class of rooted trees in
%$T(\mathbf{m})$ choose a  tree to represent this class. 
%For each pair of planar trees $s$ and $t$ that are identical
%as labeled rooted trees, for each collection $C$, there is a canonical
%isomorphism of rooted trees from $s$ to $t$. This induces an
%isomorphism of $\Aut(\mathbf{m})$-modules from
%$\bigotimes_{v\in\vert(s)}C(\leg(v))$ to
%$\bigotimes_{v\in\vert(t)}C(\leg(v))$.

Let $C$ be a collection. We define an inclusion of vector spaces from
the tensor coalgebra $T'(\pmb{\oplus}C)$ on $\pmb{\oplus}C$ into a big space  
\begin{equation}
T'(\pmb{\oplus}C) \overset{i}{\longrightarrow} \bigoplus_{t \text{ planar}}
\bigotimes_{v\in\vert(t)} C(\leg(t)).\label{eq:idef}
\end{equation}
The map $i$ on the summand $C(\mathbf{n}_1)\otimes\ldots\otimes
C(\mathbf{n}_m)\subset T'(\pmb{\oplus}C)$ has a value in the summand
corresponding to a planar tree $t$ (at the right hand side) for each
planar tree $t$. The image in the summand of $t$ is non-zero only
if $\leg(v_i) = \mathbf{n}_i$ for all $i$. If $i$ is non-zero, it is
given by the identification of $C(\mathbf{n}_1)\otimes\ldots\otimes
C(\mathbf{n}_m)$ with $C(\leg(v_1))\otimes\ldots\otimes C(\leg(v_m))$,
where the identification $C(\mathbf{n}_i)\longrightarrow C(\leg(v_i))$ is
the planar structure on $t$ (i.e. the isomorphism $\mathbf{n}_i
\longrightarrow \leg(v_i)$).

\subsection{Symmetrization}

We can identify the symmetric coalgebra $S'(\pmb{\oplus}C)$ with
the invariant part  of the tensor coalgebra  $T'(\pmb{\oplus}C)$ (with
respect to permutation of tensor factors). 
Note that this does not involve the action of tree automorphisms.
The restriction of the map $i$ in Formula (\ref{eq:idef}) to
$S'(\pmb{\oplus}C)$ is denoted $i$ as well. 

\begin{Pp}\label{Lm:mapStoT}
The map $i$ induces a map (again denoted by $i$) from
$S'(\pmb{\oplus}C)$ into the total space $\pmb{\oplus}(T'C)$ of the cofree
pseudo cooperad on $C$.  That is, $i$ defines a map 
\[
S'(\pmb{\oplus}C) \overset{i}{\longrightarrow} \pmb{\oplus}(\bar{T}'C)
\]
into the total space of $T'C$.
\end{Pp}
\begin{Pf}
Forgetting the planar structure defines a map of $\Aut(\leg(t))$-modules
\begin{equation}\label{eq:iinvpf}
S'(\pmb{\oplus}C) \overset{i}{\longrightarrow}
\pmb{\oplus}\left(
\bigoplus_{t \text{planar}}
\bigotimes_{v\in\vert(t)}C(\leg(v_i))
\right)
\overset{\psi}{\longrightarrow} 
\pmb{\oplus}\left(\bigoplus_t\bigotimes_{v\in\vert(t)}
C(\leg(v_i))\right).
\end{equation}
It now suffices to show that the image 
of the composition in Equation (\ref{eq:iinvpf}) is
invariant under the $\Aut(t)$-action for each rooted tree $t$.
We study $\psi\circ i$ on  $sum_{\sigma\in S_m}(p_{\sigma1},\ldots, p_{\sigma
m}) \in S^m(\pmb{\oplus}C)$ with $p_j\in C(n_j)$. The map
(\ref{eq:iinvpf}) satisfies 
%\begin{equation}\label{eq:iimpf}
%\sum_{\sigma\in S_m}(p_{\sigma1},\ldots, p_{\sigma m}) 
%\longmapsto \sum_{j=1}^m \sum_{(s_i)_{i=1}^{n_j}} \sum_{(P_i)_{i=1}^{n_j}}
%(s_1(P_1)\circ_1(\ldots(s_{n_j}(P_{n_j})\circ_{n_j}v_0(p_j))...)),
%\end{equation}
\begin{equation}\label{eq:iimpf}
\sum_{\sigma\in S_m}(p_{\sigma1},\ldots, p_{\sigma m}) 
\longmapsto \sum_{j=1}^m \sum_{(s_i)_{i=1}^{n_j}} \sum_{(P_i)_{i=1}^{n_j}}
((\cdots (v(P_j)\circ_{n_j} s_{n_j}(P_{n_j}))\ldots )\circ_1 s_1(P_1)),
\end{equation}
where the second sum is over $n_j$-tuples
of planar trees 
$(s_1,\ldots, s_{n_j})$ such that $\sum_k|\vert(s_k)|=m-1$ (for the
sake of simplicity of the formula a planar tree might be empty here) 
and the third sum is over ordered partitions $(P_1,\ldots,P_{n_j})$ of 
$\{1,\ldots,m\}-\{j\}$ such that $|P_k| = |\vert(s_k)|$ for all
$k$. The symbolic notation  
$((\cdots (v(P_j)\circ_{n_j} s_{n_j}(P_{n_j}))\ldots )\circ_1 s_1(P_1))$
denotes the sum of all terms where the vertices of $s_k$ are labeled
labeled by $p_l$ with $l\in P_k$, and the root vertex $v_0$
labeled by $p_j$. The summand  of $(s_1,\ldots,s_{n_j})$ corresponds
to the tree with root $r$ such that $\leg(r)= \{0,\ldots, n_j\}$ and
with planar tree $s_k$ grafted to the root along leg $k$.

We prove the invariance by induction on the number of
vertices of a tree.
Let $n\in \NN$ and suppose that the image of (\ref{eq:iinvpf}) is
contained in the invariants for each rooted tree with $<n$
vertices.
Let $s$ be a rooted tree such that  $|\vert(s)| = n$.
One can write $\Aut(s)$ as a crossed product of
$\bigoplus_i\Aut(s_i)$ with the subgroup of $\Aut(\leg(v_0))$ that
permutes those legs $l_i$ such that the corresponding $s_i$ are isomorphic.
Consequently, it suffices (by the induction hypothesis) to show that
the image is invariant with respect to elements $\hat\tau\in\Aut(s)$
defined by a permutation $\tau\in \Aut(\leg(v_0))$.
This will be done for each summand $j$ in Formula (\ref{eq:iimpf}).
Let $\tau\in \Aut(\leg(p_j))$. Then, $\tau$ applied to a summand of
the right hand side of Formula (\ref{eq:iimpf}) gives
\begin{equation}\label{eq:ieqfor}
%(s_{\tau^{-1}1}(P_{\tau^{-1}1})\circ_1(\ldots
%(s_{\tau^{-1}n_j}(P_{\tau^{-1}n_j})\circ_{n_j}v_0(p_j))...)),
((\cdots (v(P_j)\circ_{n_j} s_{\tau^{-1}n_j}(P_{\tau^{-1}n_j}))\ldots
  )\circ_1 s_{\tau^{-1}1}(P_{\tau^{-1}1})),
\end{equation}
where we compensate for the omitted action of $\tau$ on $p_j$ by 
not permuting the legs of the tree according to the induced permutation.  
This gives a permutation of summands and equivariance follows.
\end{Pf}
\begin{Pp}\label{Pp:eqmapStoT}
The map $i$ of Proposition \ref{Lm:mapStoT} descends to the coinvariants
with respect
to the action on the collection $C$. That is, $i$ defines a map 
\[
S'(\pmb{\oplus}_SC) \overset{i}{\longrightarrow} \pmb{\oplus}_S(\bar{T}'C).
\]
\end{Pp}
\begin{Pf}
This follows by the equivariance property of cocompositions in $\bar
T'C$, which we can express by compensating for the action of $\tau$ on
the label of a vertex of a tree $t$ by a permutation of $\leg(t)$,
as in Formula (\ref{eq:ieqfor}).  
\end{Pf}

%
% L infty Algebras
%

\subsection{$L_{\infty}$-Algebras}

An \note{$L_\infty$-algebra} (Lada and Stasheff \cite{LaSta:Linf},
et. al.) is a dg vector space $\Gg$, together with a
differential $\pd$ (i.e. a coderivation s.t. $\pd^2=0$) on the free
cocommutative coalgebra $S'(\Gg[-1])$. 
The differential is completely determined by its restrictions
$l_n:S^n(\Gg[-1])\longrightarrow \Gg[-1]$,
where $S^n(\Gg[-1])$ of course denotes the elements of tensor degree
$n$. A note{morphism of $L_\infty$-algebras} (or $L_\infty$-map)
$f:\Gg \longrightarrow \Gh$ is a morphism of the cofree coalgebras
compatible with the differentials. It is determined by its restrictions
$f_n:S^n(\Gg[-1])\longrightarrow \Gh[-1]$. An $L_\infty$-map is a
\note{quasi isomorphism} if $f_1:\Gg[-1]\longrightarrow \Gh[-1]$ induces an
isomorphism in cohomology.

\begin{Lm}\label{Lm:STcoder}
Let $C$ be a collection and let $D$ be a coderivation on
$T'C$. Then the diagram
\[
\xymatrix{
S'(\pmb{\oplus}C) \ar[d]_D \ar[r]^i& \pmb{\oplus}(T'C) \ar[d]^{\oplus D}\\
S'(\pmb{\oplus}C) \ar[r]^i& \pmb{\oplus}(T'C)
}
\]
commutes, where $D:S'\pmb{\oplus}C\longrightarrow S'\pmb{\oplus}C$ is the
coderivation defined by the composition
\[
\xymatrix{
S'(\pmb{\oplus}C)\ar@{.>}[d] \ar[r]^i 
& \pmb{\oplus}(T'C) \ar[d]^{\pmb{\oplus}D} \\ 
\pmb{\oplus} C 
& \pmb{\oplus}(T'C) \ar@{->>}[l].
}
\]
\end{Lm}
\begin{Pf}
This is a direct corollary of the definition of coderivation for coalgebras
and cooperads.
\end{Pf}

\begin{Tm}\label{Tm:Linfty}
Let $P$ be an operad up to homotopy. Then
\begin{enumerate}
\item there exists a natural $L_\infty$-algebra structure on
$\pmb{\oplus} P$. 
\item 
this structure descends to the quotient $\pmb{\oplus}_{S} P$. 
\end{enumerate}
\end{Tm}
\begin{Pf}
Observe that the differential $\pd$ on $T'(P[-1])$ induces a derivation $d$ on
$S'(\pmb{\oplus}P[-1])$ according to Lemma \ref{Lm:STcoder}.
Regarding part \textsl{(i)}, we need to show that $d^2=0$. 
Since $d$ and $i$ commute, we have a
commutative diagram
\[
\xymatrix{
S'\left(\pmb{\oplus} P[-1]\right)\ar@/^2pc/@{.>}[rrdd]\ar[r]^i\ar[d]^d &
\pmb{\oplus}(T'(P[-1])) \ar[d]^d&\\
S'\left(\pmb{\oplus} P[-1]\right)\ar[r]^i\ar[d]^d & \pmb{\oplus}(T'(P[-1]))
\ar[d]^d&\\
S'\left(\pmb{\oplus} P[-1]\right)\ar[r]^i & \pmb{\oplus}(T'(P[-1]))
\ar@{->>}[r]^\pi &\pmb{\oplus}P[-1]}
\]
with the composition from $S'\left(\pmb{\oplus} P[-1]\right)$ at the
top left to $\pmb{\oplus}P[-1]$ (dotted arrow) equal to
$0$. Therefore, $d^2 = 0$ on
$S'(\pmb{\oplus}P[-1])$. This shows \textsl{(i)}.

Since the differential on $T'(P[-1])$ preserves the symmetric group
action (cf. $\circ_t$ operations),
it induces a differential on $\pmb{\oplus}_ST'(P[-1])$. The map $i$ is
well defined as a map $i:S'(\pmb{\oplus}_SP[-1])\longrightarrow
\pmb{\oplus}_ST'(P[-1])$. The induced differential on $
S'(\pmb{\oplus}_SP[-1]$  commutes with the map $i$. Part \textsl{(ii)}
follows
by the commutation relation of the $\circ_t$ with the
$Aut(\leg(t))$-actions. 
\end{Pf} 
\begin{Rm}
It is instructive to write the higher operations $l_n$ in the
$L_\infty$-algebra
structure of the previous theorem. 
\begin{equation}
l_n = \sum_{\tiny\begin{array}{c} t \text{ planar},\\
|\vert(t)|=n\end{array}} 
(\circ_t)\circ i: S^n \left(\pmb{\oplus}P[-1]\right)
\longrightarrow \pmb{\oplus}P[-1],
\end{equation}
The signs are induced by the Koszul convention for the
$|\vert(t)|$-ary operations $\circ_t$.
\end{Rm}
%
% A infty
%
\subsection{Examples: $A_\infty$-Algebras and Operads}

This paragraph treats some special cases of the $L_\infty$-algebras
constructed in the previous section.
Firstly we restrict our attention to $A_\infty$-algebras, where we obtain a
well known result as a corollary to Theorem \ref{Tm:Linfty}. This
example is based on the observation that $A_\infty$-algebras
are exactly operads up to homotopy concentrated in arity 1.
Secondly, we recover the Lie algebra structure on the total space
$\pmb{\oplus}P$ of an operad $P$.  

An \note{$A_\infty$-algebra} (Stasheff \cite{Sta:Hspace}) is a vector
space $A$, together with a 
differential $\pd$ (i.e. a coderivation s.t. $\pd^2=0$) on the free
 coalgebra $T'(A[-1])$. 
The differential is completely determined by its restrictions
\[
m_n:T^n(A[-1])\longrightarrow A[-1],
\]
where $T^n(A[-1])$ of course denotes the elements of tensor degree
$n$.
\begin{Lm}
Let $(T'(P[1]),d)$ be an operad up to homotopy. Then the restriction of $d$
to the part of arity 1 makes $P(1)$ an $A_\infty$-algebra.
\end{Lm}
\begin{Pf}
Let $A = P(1)$. the differential on the tree cooperads $T'(P[-1])$
restricts to a 
differential on the subspace of $T'(P[-1])(1)$ consisting of the
summands corresponding to trees in which $|\leg(v)|=2$ for each vertex
$v$. This sub-cooperad equals the cofree coalgebra
$T'(A[-1])$ on $A[-1]$. Moreover, the differential is a coderivation
since it is a coderivation on $T'(P[-1])$. This proves the result.
\end{Pf}

\begin{Ex}[Lada and Markl \cite{LaMa:Linf}]
When we start out from $P$ concentrated in arity 1, the previous
results shows that each $A_\infty$-algebra becomes an $L_\infty$-algebra 
with what one might call higher order commutators. This
generalizes the construction of the commutator Lie algebra with the bracket
$[a,b]=ab-ba$ of an associative algebra $A$. 
%For convenience of the
%reader we  sketch the proof of this special case of Theorem
%\ref{Tm:Linfty} in terms of the classical structures.
%One can construct the symmetric coalgebra $S'(A[-1])$ as the invariant part
%(w.r.t. permutation of tensor factors) of
%the tensor coalgebra $T'(A[-1])$ in each tensor degree. Since the
%differential $d$ is a coderivation it is 
%determined by its projection $T'(A[-1])\longrightarrow A$. Composition
%with the inclusion $S'(A[-1])$ shows that it restricts to a
%coderivation of $S'(A[-1])$ (which is cogenerated by $A$).  
%
%To see the connection with Theorem \ref{Tm:Linfty}, observe that
One can
write the operations $l_n$ as the symmetrization of $m_n$. 
%Since all operations are of arity one both $L_\infty$-structures coincide.
\end{Ex}

The first part of the following proposition has been around for a long
time. 
It suffices to use an argument of Gerstenhaber \cite{Gerst:Ring} (predating
the very 
definition of operads!). The second part first appeared in Kapranov and
Manin \cite{KaMa:ModMor}, though for example Balavoine \cite{Bal:Def} 
already gives a sufficient computation. But he failed to state a
result in this generality. 
\begin{Cr}\label{Pp:KapMan}
Let $P$ be an operad. Then
\begin{enumerate}
\item The total space $\pmb{\oplus} P$ is a Lie algebra with
respect to the commutator of the (non-associative) multiplication
$p \circ q = \sum_{i=1}^n p \circ_i q$ for $p\in P(n)$ and $q\in P(m)$. 
\item The Lie algebra structure descends to  quotient $\pmb{\oplus}_SP$.
\end{enumerate}
\end{Cr}

\subsection{Homotopy Homomorphisms}

This section studies morphisms of (co)operads up to homotopy. We are
interested in the interaction between 
these maps and the $L_\infty$-structures described in the previous
paragraph. 

\begin{Df}
A \note{morphism of operads up to homotopy} $\phi:P\longrightarrow Q$  is a
morphism of dg cooperads $\phi:T'(P[-1],d)\longrightarrow
T'(Q[-1],d)$.

Dually, a \note{morphism  of cooperads up to homotopy}
$\psi:A\longrightarrow C$ is a morphism of dg operads
$\psi:T(A[1],d)\longrightarrow T(C[1],d)$.
\end{Df}

The sequel of this section focuses on the operad up to homotopy case
and leaves the
formulation of the dual statements as an exercise. The adjuction of
$T'$ and the forgetful functor assures that $\phi$ is determined by
maps $\phi(t):(P[-1])(t)\longrightarrow P[-1]$. The condition that 
$\phi$ is compatible with the differential can be described in terms of
conditions about compatibility with the $\circ_t$ operations:
\begin{equation}
\sum_{s\subset t} \pm \phi(t/s) \circ (\circ_s) -
(\circ_{t/s})\circ\phi(s) = 0,
\end{equation}
again with the signs depending on the trees as in Formula
(\ref{eq:pdsquare}) .

\begin{Lm}\label{Lm:STmorph}
Let $C$ and $E$ be collections, and let 
$\phi:T'C\longrightarrow T'E$ be a morphism of cooperads. Then the diagram
\[
\xymatrix{
S'\pmb{\oplus}C \ar[d]_{\tilde{\phi}} \ar[r]^i& \pmb{\oplus}T'C
\ar[d]^{\pmb{\oplus}\phi}\\
S'\pmb{\oplus}E \ar[r]^i& \pmb{\oplus}T'E
}
\]
is commutative, where $\tilde{\phi}:S'\pmb{\oplus}C\longrightarrow
S'\pmb{\oplus}E$ is the coalgebra morphism induced by
\[
\xymatrix {
S'(\pmb{\oplus}C) \ar[r]^i\ar@{.>}[d] & {\pmb{\oplus}(T'C)} 
\ar[d]^{\pmb{\oplus}\phi} \\
(\pmb{\oplus}E) & {\pmb{\oplus}(T'E)} \ar@{->>}[l].
}
\]
\end{Lm}
\begin{Pf}
This is clear from the definition of coalgebra and cooperad morphism. 
\end{Pf}

\begin{Tm}\label{Tm:OpdLinfmap}
Let $\phi:P\longrightarrow Q$ be a morphism of operads up to
homotopy. Then $\phi$ induces morphisms of $L_\infty$-algebras
\[
\pmb{\oplus}P \longrightarrow \pmb{\oplus}Q
\qquad \text{ and }\qquad
\pmb{\oplus}_SP \longrightarrow \pmb{\oplus}_SQ.
\]
\end{Tm}
\begin{Pf}
The universal property of $S'(\pmb{\oplus}Q[-1])$ defines a coalgebra
morphism $\phi$ as in Lemma \ref{Lm:STmorph}. 
It suffices to check that the coalgebra map $\phi$ commutes with the
differential. Due to the coderivation property (Lemma
\ref{Lm:STcoder}) it suffices to show that $\pi\circ d \circ \phi =
\pi \circ \phi \circ d$, where 
$\pi:S'(\pmb{\oplus} Q[-1]) \longrightarrow \pmb{\oplus}Q$ is the
natural projection. The result now
follows since $i$ and $\phi$ commute by Lemma \ref{Lm:STmorph} and
since $\phi$ is equivariant with respect to the $S_n$-actions. 
\end{Pf}

%
%   M O D E L S
%

\section{Auxiliary Results}

In order to proceed, we need to recall some results.  Most
results in this section are well known. Nevertheless we feel inclined
to spell out some details on $L_\infty$-algebras and the Maurer-Cartan
equation.

\subsection{Model Category of Operads}

We recall the model category (cf. Quillen \cite{Quil:Htpy} and Hovey
\cite{Hov:Model}) of operads first defined in Hinich
\cite{Hin:Htpy}.  The weak equivalences are quasi isomorphisms,
fibrations are surjections, and the generating cofibrations are
extensions by free operads with a differential.

We recall some facts on minimal models of 1-reduced operads, following
Markl \cite{Mar:Model}.
Let $P$ be a 1-reduced operad. A \note{model} $A$ of $P$ is a cooperad up to
homotopy $A$, together with a quasi isomorphism $T(A[-1],\pd)\longrightarrow
P$.  A model $M= (T(A[-1]),\pd)$ of $P$ is a \note{minimal model} if the
differential $\pd$ is
minimal in the sense that the induced internal differential on $A$
vanishes. 
A model $M= (T(A[-1]),\pd)$ of $P$ is a \note{strict model} if $A$ is
strict. That is, the differential $\pd$ defines a cooperad structure on $A$.

\begin{Tm}[Markl\cite{Mar:Model}]
Every 1-reduced operad has a minimal model. Moreover, the minimal
model is unique up to isomorphism of cooperads up to homotopy.
\end{Tm}

\begin{Rm}[Hinich \cite{Hin:Htpy}] 
Operads in $k\dgVect$ form a model category. 
Cofibrant objects in the model category are exactly retracts of operads
associated to homotopy cooperads. Models and minimal models are cofibrant
objects (of a special kind!) in this model category. To see this, use
that fact that the generating cofibrations are free operads on a cofibrant
collection, together with a differential. Since the characteristic of $k$ is
0, every collection is cofibrant. A minimal model $(T(A[-1]),\pd)$ of an
operad $P$
is a cofibrant replacement for $P$ if the map $TA[1]\longrightarrow P$ is a
fibration (i.e. is
surjective). It suffices to ask the internal differential of $P$ to be 0. 
\end{Rm}

%
% Maurer-Cartan
%

\subsection{Maurer-Cartan Equation}

Throughout this section, 
let $(S'(\Gg[-1]),\pd)$ be an $L_\infty$-algebra and let $\phi\in\Gg^1$
be a solution of the \note{(generalized) Maurer-Cartan equation} 
\begin{equation}\label{eq:MC}
\sum_{n\geq 1}l_n(\phi^{\otimes n}) = 0.
\end{equation}
Note that in the case of a Lie algebra this equation takes the form 
\[ 
d\phi + \frac{1}{2}[\phi,\phi] = 0,
\]
where the factor $\frac{1}{2}$ comes from using coinvariants instead
of invariants. We will need some lemmas on these structures. The first
lemma shows how a solution of the Maurer-Cartan equation can be
interpreted as a perturbation of the $L_\infty$-structure. 
The second lemma describes how a solution of the Maurer-Cartan equation
can be transported along an $L_\infty$-map. The third lemma shows that
in this situation there exists a morphism of the perturbed
$L_\infty$-algebras. These results should be well known to the experts
in the field. Our formulas differ a bit from the usual ones due to
consequent identification of $S'(\Gg[-1])$ with invariants (as opposed
to coinvariants). 

\begin{Cv}
In this section we work formally. That is, the formulas
contain a priori infinite series and we do not worry about
convergence. All results are understood to hold modulo convergence of
the relevant series. A sufficient condition for convergence will be stated
in the next section.
\end{Cv}

First we elaborate a bit on the exact definition of the operations in
$L_\infty$-algebras.
The differential $\pd$  on $(S'(\Gg[-1]),\pd)$ can be described in terms of
the maps
$l_n^k:S^{n+k} \longrightarrow S^{k+1}$, where 
\[
l_n^k = \text{Sh}_{1,k}\circ(l_n\otimes \id^{\otimes k}).
\]
Here $\text{Sh}_{p,q}$ denotes the sum of all $(p,q)$-shuffles
permuting tensor factors. These serve to assure that the result is
again invariant under the action of the symmetric group. 
The square-zero equation reads
\begin{equation}
\label{eq:Linfdiff}
\sum_{n+k = m} l_{k+1}\circ l_n^k = 0 \qquad \text{ for all } m.
\end{equation}
An $L_\infty$-map $f$ satisfies  
\begin{equation}
\label{eq:Linfmap}
\sum_m\sum_{(n_1,\ldots,n_m)} l_m\circ(f_{n_1}\otimes\ldots\otimes
f_{n_m}) = \sum_{k+ p=n} f_{k+1}\circ l^k_{p},
\end{equation}
where the sum over $(n_1,\ldots,n_m)$ assumes $n_1+\ldots+ n_m=n$.
The only difficulty in establishing the results below is in careful
bookkeeping with respect to these maps. To do this in a satisfactory
way we need some properties of the shuffle maps. The shuffles product
is associative and graded commutative in the sense that
\begin{align*}
\Sh_{p+q,r}\circ(\Sh_{p,q}\otimes\id^{\otimes r}) &= \Sh_{p,q+r}\circ
(\id^{\otimes p}\otimes \Sh_{q,r})\\
\Sh_{p,q} &= \Sh_{q,p} \circ \tau,
\end{align*}
where $\tau:V^{\otimes p}\otimes V^{\otimes q} \longrightarrow
V^{\otimes q}\otimes V^{\otimes p}$ is the twist of tensor factors.
Note that the usual signs are hidden is the action of the shuffles.

The lemmas below are obtained by applying the definitions above and
the properties of the shuffle product.
\begin{Lm}\label{Lm:perturbLinf}
Let  $(S'(\Gg[-1]),\pd)$ be a $L_\infty$-algebra and
let $\phi\in\Gg^1$ be a solution of the
Maurer-Cartan equation.  Then
$(S'(\Gg[-1]),\pd_\phi)$ is an $L_\infty$-algebra with the
differential given by 
\[
\tilde{l}_j(x_1,\ldots,x_j) = \sum_{p\geq 0}
l_{j+p}(\text{Sh}_{p,j}(\phi^{\otimes
p};x_1,\ldots,x_j).
\] 
\end{Lm}

\begin{Lm}\label{Lm:moveMC}
Let $f:(S'(\Gg[-1]),\pd)\longrightarrow (S'(\Gh[-1]),\pd)$ be a map of
$L_\infty$-algebras and let $\phi$ be a solution of the Maurer-Cartan
equation in $(S'(\Gg[-1]),\pd)$. Then 
\[
\psi = \sum_n f_n(\phi^{\otimes n})
\]
is a solution of the Maurer-Cartan equation in $(S'(\Gh[-1]),\pd)$.
\end{Lm}

\begin{Lm}\label{Lm:perturbmap}
Let $f:(S'(\Gg[-1]),\pd)\longrightarrow (S'(\Gh[-1]),\pd)$ be a map of
$L_\infty$-algebras. Let $\phi$ be a solution of the Maurer-Cartan
equation in $(S'(\Gg[-1]),\pd)$ and let $\psi$ be the induced solution
of the Maurer-Cartan equation on $(S'(\Gh[-1]),\pd)$. Then  the formula
\[
\tilde{f}_n(x_1,\ldots,x_n) =\sum_{p\geq 0}
f_{n+p}(\text{Sh}_{p,n}(\phi^{\otimes p};x_1,\ldots,x_n)).
\]
defines a morphism of $L_\infty$-algebras
$\tilde f = f_\phi:(S'(\Gg[-1]),\pd_\phi)\longrightarrow
(S'(\Gh[-1]),\pd_\psi)$.
\end{Lm}

\subsection{Convergence and Quasi Isomorphisms}

The lemmas above need some extension. We study
how the constructions above behave under quasi
isomorphism. In order for this to work we need some mild assumptions. 

For an $L_\infty$-algebra $(S'(\Gg[-1]),\pd)$ we denote by $d$ the
internal differential induced by $\pd$. If $\phi$ is a solution of the
Maurer-Cartan equation we denote by $\tilde{d}$  
the internal differential of the $\phi$-perturbed $L_\infty$-algebra.

An $L_\infty$-algebra is called \note{bigraded} if we can write $\Gg =
\bigoplus_{n\in\NN}\Gg_n$ as dg vector spaces, and the differential
$\pd$ preserves the induced grading on the symmetric coalgebra. If
$\Gg$ is a bigraded $L_\infty$-algebra, we can construct the
completion of $\Gg$ in the usual way. The completion
$\hat\Gg=\hat{\bigoplus}_n\Gg_n = {\prod}_n\Gg_n$ is an $L_\infty$-algebra
with respect to the completed tensor product. We call completions of
bigraded $L_\infty$-algebras \note{bigraded complete}
$L_\infty$-algebras. 
A morphism of bigraded complete $L_\infty$-algebras is assumed to
preserve the second grading. 
\begin{Rm}
The completeness in the definition above 
assures that all the series in the previous section are well defined for
bigraded complete $L_\infty$-algebras with a solution $\phi$ of the
Maurer-Cartan equation in $\hat{\bigoplus}_{m>0}\Gg_m$.
\end{Rm}
\begin{Pp}
Let $(S'(\Gg[-1]),\pd)$ and $(S'(\Gh[-1]),\pd)$ be bigraded
$L_\infty$-algebras. Let $\phi\in \hat{\bigoplus}_{m>0}\Gg_m$  be a solution to
the Maurer-Cartan equation in the completed $L_\infty$ algebra. Let
$f:(\Gg,d)\longrightarrow (\Gh,d)$ be a quasi
isomorphism of bigraded $L_\infty$-algebras. Then the map
$f_\phi$ of Lemma \ref{Lm:perturbmap} is a quasi isomorphism of the
completed $L_\infty$-algebras.
\end{Pp}
\begin{Pf}
This is a direct corollary of the two lemmas below, since the
filtration $F^{\geq n}\hat\Gg$ defined in the proof of Lemma
\ref{Lm:specseq} is exhaustive, weakly convergent and complete
(cf. McCleary \cite{McCl:Spec}).
\end{Pf}
\begin{Lm}\label{Lm:specseq}
Let $(S'(\Gg[-1],\pd)$ be a bigraded $L_\infty$-algebra  and
let $\phi\in \hat{\bigoplus}_{m>0}\Gg_m$ be a solution of the
Maurer-Cartan equation in the completed $L_\infty$-algebra. 
Then there exists a spectral sequence
\[
E_1^{pq}(\hat\Gg) = \left(F^pH(\hat\Gg,d)/F^{p+1}H(\hat\Gg,d)\right)^{p+q},
\]
that converges to $H(\hat\Gg,\tilde{d})$. 
\end{Lm}
\begin{Pf}
The spectral sequence is the spectral sequence defined by the
filtration of $(\hat\Gg,\tilde{d})$ given by 
\[
F^{\geq n}\hat\Gg = \hat{\bigoplus}_{m\geq n}\Gg_m. 
\]
The term of $\tilde{d}(x)= \sum_k l_{k+1}(\text{Sh}_{k,1}(\phi^{\otimes
k}; x))$ that involves $\phi^{\otimes k}$ increases the grading 
and thus contributes to the differential only in $E_r$, where 
$r\geq k$. Then $E_0^{pq} = (F^p\hat\Gg/F^{p+1}\hat\Gg)^{p+q}$, and $d_0$ is
the internal differential $d$ of $\hat\Gg$. Since $\Gg$ is bigraded,
the result follows.
\end{Pf}

\begin{Lm}\label{Pp:specseq}
Let $f:(\Gg,d)\longrightarrow (\Gh,d)$ be a quasi
isomorphism of bigraded $L_\infty$-algebras. 
Let $\phi\in \hat{\bigoplus}_{m>0}\Gg_m$  be a solution to
the Maurer-Cartan equation and let $\psi = \sum_nf_n(\phi^{\otimes n})$.
Then the morphism of $L_\infty$-algebras
$\tilde f:\hat\Gg_\phi\longrightarrow\hat\Gh_\psi$ induces a morphism
of spectral sequences $E_k^{**}(\hat\Gg_\phi) \longrightarrow
E_k^{**}(\hat\Gh_\psi)$, which is an isomorphism for $k\geq 1$. 
\end{Lm}
\begin{Pf}
Since the $f$ is a morphism of bigraded $L_\infty$-algebras, the
morphism $\tilde f_1:\Gg_\phi\longrightarrow \Gh_\psi$ commutes with the
internal differentials $\tilde{d}$ in $\hat\Gg_\phi$ and $\hat\Gh_\psi$, and
preserves the filtration. Therefore $\tilde f_1$defines a morphism of
spectral sequences. The part of $\tilde f$ which has degree 0 with
respect to the second grading is $f_1$ since
$\phi\in\hat{\bigoplus}_{m\geq 1}\Gg_m$. Thus the map on $E_0$ is given
by $f_1$. Since $f_1$ is a quasi isomorphism, the explicit form of
$E_1$ given in Lemma \ref{Lm:specseq} shows that the map on $E_1$ is 
\[
E_1(\hat\Gg) =
H(\hat\Gg,d)\overset{H(f_1)}{\longrightarrow}H(\hat\Gh,d) =
E_1(\hat\Gh).
\]
This suffices to prove the result (cf. McCleary \cite{McCl:Spec}).
\end{Pf}

\section{Cohomology and Deformations}

This section defines for any reduced operad $Q$ and each operad
$\phi:Q\longrightarrow P$ under $Q$ the operad cohomology $H_Q(P)$ of $Q$ 
with coefficients in $P$. 
We identify $H_Q(P)$ with the cotangent
cohomology defined by Markl \cite{Mar:Model} using an
$L_\infty$-structure on the cotangent complex.  We show how the
approach of Balavoine \cite{Bal:Def}, and Konstevich-Soibelman
\cite{KonSoi:Del} to deformations of algebras over an operad fit in
this framework.

\subsection{Maurer-Cartan for Convolution Operads}

The associated Lie algebra of an operad $R$ always has a second grading
given by the arity as $(\pmb{\oplus}(R))_m = R(m-1)$. 
Let $A$ be a cooperad and let $P$ be an operad. This section
examines dg operad maps from the cobar construction $B^*A$ to $P$ in
relation to the convolution operad $P^A$ and the Maurer-Cartan
equation in the completion $\hat{\pmb{\oplus}}_SP^A$,
of $\pmb{\oplus}_SP^A$.

\begin{Tm}\label{Tm:OpdMC}
Let $A$ be a cooperad and $P$ an operad.  
There is a one-one correspondence between elements
$\phi\in(\hat{\pmb{\oplus}}_SP^A)^1$ that satisfy the Maurer-Cartan
equation $\frac12[\phi,\phi] + d(\phi) = 0$ and
dg operad morphisms $\tilde{\phi}:B^*(A)\longrightarrow P$.
\end{Tm}
\begin{Pf}
Identify the collection morphisms from $A$ to $P$ of degree $n$ with
the objects of $(\hat{\pmb{\oplus}}_SP^A)^n$ (i.e. elements of internal
degree $n$) using the section $\frac{1}{n!}\sum_\sigma \sigma$.
By the universal property of the free operad, an operad morphism is
uniquely determined by a collection morphism $\tilde\phi:A[1]\longrightarrow P$ in
$k\gVect$. To complete the proof we need to see that
the condition $\frac12[\phi,\phi]+d(\phi)=0$ in $\pmb{\oplus}_SP^A$ is
equivalent to the induced map $\tilde\phi:T(A[1]) \longrightarrow P$ being
compatible with the differentials. Of course, here $d(\phi)$ denotes the
differential induced by the convolution operad $P^A$ (cf. Lemma
\ref{Lm:convol}). 

Denote the total differential on $B^*A$ induced by the cooperad
structure by $\pd_A$ and the internal differential of $A$ by $d$. Let
$\pd = \pd_A + d$ be the total differential. 
Compatibility of $\tilde\phi$ with the differential can be stated as 
\[
\tilde\phi\circ\pd_A = d\circ \tilde\phi - \tilde\phi\circ d,
\]
where the internal differential of $P$ is denoted $d$ as well. 
The restriction of $\tilde\phi\circ\pd_A$ to $A$ is by definition of the
Lie bracket $\fr12 [\phi,\phi]$. 
\end{Pf}
\begin{Rm}
As Martin Markl observed, Theorem \ref{Tm:OpdMC} gives a Lie algebraic
framework for the twisting functions of Getzler and Jones \cite{GetzJon:Opd}.
\end{Rm}
We can deduce an invariance property of these convolution operads from
generalities on model categories and our previous results.
\begin{Pp}\label{Pp:quasiPAPC}
Let $A$ and $C$ be cooperads.
Suppose that $B^*A$ and $B^*C$ are cofibrant replacements for the same
operad $Q$, then there exists a pair of quasi isomorphisms of
$L_\infty$-algebras
\[
\begin{split}
&\hat{\pmb{\oplus}}_SP^A \longrightarrow \hat{\pmb{\oplus}}_SP^C \qquad \text{ and }\\
&\hat{\pmb{\oplus}}_SP^C \longrightarrow \hat{\pmb{\oplus}}_SP^A
\end{split}
\]
that are inverse in cohomology.
\end{Pp}
\begin{Pf}
We first construct
$\eta:B^*A\longrightarrow B^*C$ and $\zeta:B^*C\longrightarrow
B^*A$ of dg operads such that $H\eta$ and $H\zeta$ are each others inverse. 
Both models are cofibrant in the model category of operads. In the solid
square
\[
\xymatrix{ 0 \ar[d]\ar[r] & B^\ast A \ar[d]\ar@{.>}@<.5ex>[dl]\\
B^\ast C\ar[r]\ar@{.>}@<.5ex>[ur] & Q
}
\]
the top and left arrows are cofibrations and the bottom and right
arrows are trivial fibrations. Consequently,
the two dotted diagonals $\eta$ and $\zeta$ exist.

The maps $\eta$ and $\zeta$ induce  homotopy homomorphisms 
$\eta: P^C \longrightarrow P^A$ and
$\zeta: P^A \longrightarrow P^C$.
These are quasi isomorphisms by the K\"unneth formula. 
The induced quasi isomorphisms of $L_\infty$-algebras (cf. Theorem
\ref{Tm:OpdLinfmap}) are the desired maps.
\end{Pf}

\subsection{Operad Cohomology}

This section introduces a cohomology theory that we will show to
contain the approaches of Markl \cite{Mar:Model}, Balavoine
\cite{Bal:Def}, and Kontsevich and 
Soibelman \cite{KonSoi:Del} on the deformation of algebras 
over operads. The present construction is done purely at the operad
level and thus generalizes to a dg Lie algebra controlling
deformations of operad maps.

A straightforward application of Lemma \ref{Lm:perturbLinf} in the case of
convolution operads as studied in the previous section yields the
following:
\begin{Cr}
Let $A$ be a cooperad and let $P$ be an operad.
Let $\hat\phi:B^*A\longrightarrow P$ be a map of dg operads. Then 
$D = d + [\phi,-]$ 
makes $L_A(P) = (\hat{\pmb{\oplus}}_SP^A,D)$ a dg Lie algebra, where the
correspondence between $\phi$ and $\hat\phi$ is the correspondence of
Theorem \ref{Tm:OpdMC}.
\end{Cr}

\begin{Df}\label{Df:Opcohom}
An operad $Q$ is \note{1-reduced} if $Q(0) = 0$ and $Q(1)$ is a 1
dimensional space spanned by the identity.
Let $A$ be a cooperad, let $Q$ and $P$ be operads, such that $B^*A$ is
a cofibrant resolution of $Q$ and $Q$ is reduced. Let
$\phi:Q\longrightarrow P$ be a map of operads.  The cohomology of $Q$ with
coefficients in $P$ is the cohomology of the dg Lie algebra
$L_A(P)$. The cohomology of $L_A(P)$ is
called the \note{operad cohomology} of $Q$ with coefficients in $P$
and is denoted $H^*_Q(P)$.
\end{Df}

We now show that the cohomology $H^*_Q(P)$ is independent of
the cofibrant resolution $A$. This justifies the notation that
suppresses $A$. 

\begin{Tm}
Let $Q$ be a 1-reduced operad, and let $\hat{\phi}:Q\longrightarrow P$ be an
operad under $Q$. Any two complexes of the form $L_A(P)$ as in
\ref{Df:Opcohom} are quasi isomorphic. Consequently, the operad
cohomology $H_Q(P)$ is independent of the choice of a strict model $A$.  
\end{Tm}
\begin{Pf}
Let $\pi_A:B^*A\longrightarrow Q$ and $\pi_C:B^*C\longrightarrow Q$ be
two strict cofibrant replacements for $Q$.  
Consider the convolution operads $P^C$ and $P^A$. These are quasi
isomorphic by Proposition \ref{Pp:quasiPAPC}. Let $\phi_A$ 
(resp. $\phi_C$) be the solution of the Maurer-Cartan equation in the
completed $\pmb{\oplus}_SP^A$ (resp. $\pmb{\oplus}_SP^C$) defined as in
Theorem
\ref{Tm:OpdMC} by the composition of $\hat{\phi}$ with  $\pi_A$
(resp. $\pi_C$). Since $Q$ is reduced, these solutions of the
Maurer-Cartan equation strictly increase the arity-grading. The result
now follows from Lemma \ref{Lm:perturbmap} and Proposition \ref{Pp:specseq}.
\end{Pf}

\begin{Ex}\label{Ex:recover}
It might be good to point out that we recover well
known results an special cases of the theorem above. For example, let
$Q$ be a quadratic Koszul (!) operad (cf. Ginzburg and Kapranov
\cite{GinKap:Koszul}), and let $P = \End_V$, the
endomorphism operad of a $Q$-algebra $V$. Using the strict model
$B^*A = B^*(Q^{\bot})$ in Definition \ref{Df:Opcohom}, the cohomology
$H_Q(P)$ is easily seen to be the cohomology of the $P$-algebra
$V$ (cf. Balavoine \cite{Bal:Def}, Ginzburg and Kapranov
\cite{GinKap:Koszul}).
\end{Ex}
\begin{Ex}
If $Q= \text{Ass}$ the operad of associative
algebras, and $P= \End_V$ the endomorphism operad of an associative
algebra $V$
Example \ref{Ex:recover} recovers the Hochschild cohomology complex with
coefficients in $V$. If $Q=\text{Ass}$ and $P$ is any operad with
multiplication, $H_Q(P)$ is the Hochschild cohomology of the operad $P$. 
\end{Ex}
\begin{Ex}
If $Q=\text{Lie}$ the operad of Lie algebras, and $P=\End_V$ the
endomorphism operad of a Lie algebra $V$, Example \ref{Ex:recover}
recovers the (Chevalley-Eilenberg) Lie algebra cohomology complex with
coefficients in $V$. 
\end{Ex} 

%\begin{Rm} MAKE THIS PRECISE
%Invariance under quasi isomorphism of algebras is a trivial
%consequence of the results above. Let $f: V\longrightarrow W$ be a
%quasi isomorphism of algebras. Show that we can construct a quasi
%isomorphism of $L_\infty$-algebras.
%\end{Rm}

\subsection{Cotangent Cohomology}

This section describes how we can use non-strict (and in particular
minimal) models to compute the operad cohomology $H_Q(P)$. The previous
section only shows how to 
do this if the minimal model is strict. The advantage of minimal
models is of course that they tend to be small. Moreover, the vanishing
internal differential simplifies computations. We first
construct a non-strict analogue of the convolution operad. 
The proof is the operad-analogue of the $A_\infty$-algebra structure
on the tensor product of an $A_\infty$-algebra with an associative algebra.
\begin{Lm}
Let $A$ be a cooperad up to homotopy, and
let $P$ be an operad. Then the collection $P^A$ has a natural
structure of an operad up to homotopy. 
\end{Lm} 
\begin{Pf} 
%Identify $P^A$ with the collection $P \otimes A^*$, where $*$ denotes the
%linear dual. We need to give a differential on $T(P\otimes A^*[-1])$.
%The summand corresponding to a rooted tree $t$ equals
%\begin{equation}\label{eq:tfactor}
%(P\otimes A^*)(t) =
%\left(\bigotimes_{v\in\vert(t)} P(\leg(v))\otimes
%A^*(\leg(v))[-1]\right)^{\Aut(t)}.
%\end{equation}
%We can identify $\left(\bigotimes_{v\in\vert(t)}
%A^*(\leg(v))[-1]\right)^{\Aut(t)}$ with the summand corresponding to $t$
%in $(TA[1])^*$. Here we have the differential $\pd^*$ dual to the
%defining differential of the cooperad up to homotopy $A$. 
%We define the differential of $T(P\otimes A^*)$ on the factor
%(\ref{eq:tfactor}) as
%$\gamma(t)\otimes (\circ_t)$ where $\gamma(t)$ is the operad
%composition $P(t)\longrightarrow P(\leg(t))$ and $(\circ_t)$ is the dual
%of the $\circ_t^*$ operation on $A$. Since the composition in $P$
%is associative, this is a well defined differential. Moreover, this
%differential commutes with the internal differential $d$ of $P$. The
%total differential $\pd = \pd^*+d$ makes $P^A$ an operad up to homotopy. 
%BETER BEWIJS:
Use $\circ_t(\phi_1,\ldots\phi_n) =
(\gamma_t)\circ(\phi_1,\ldots,\phi_n)\circ (\circ_t^*)$ to define
composition for trees with
at least 2 vertices. Here $\gamma_t:P(t)\longrightarrow P(\leg(t))$ is
the composition in $P$. Since the operad structure on $P$ is strict, this
defines an operad up to homotopy if we define the internal
differential as precomposition with the internal differential $d_A$ on $A$. 
If we add postcomposition with $d_P$ to the internal differential this
still defines an operad up to homotopy, since $d_P$ and composition in
$P$ commute. 
\end{Pf}

\begin{Lm}
Let $A$ be a cooperad up to homotopy, and let $P$ be an operad.
There is a 1-1 correspondence between operad maps
$\phi:(T(A[1]),\pd)\longrightarrow P$ and solutions of the
Maurer-Cartan equation in the $L_\infty$-algebra $\hat{\pmb{\oplus}}_SP^A$.
\end{Lm}
\begin{Pf}
Same argument as in the proof of Theorem \ref{Tm:OpdMC}.
\end{Pf}
Let $Q$ be a reduced operad, let $A$ be a cooperad up to homotopy such
that $\pi_A:T(A[-1]),\pd)\longrightarrow Q$ is a quasi isomorphism,
and let $\phi:Q\longrightarrow P$ be a morphism of operads. We use the
above lemma to define the complex $L_A(P)$ as the $\phi$-perturbed
$L_\infty$-algebra with respect to the perturbed internal
differential.
\begin{Tm}\label{Tm:nostrictcohom}
Let $Q$ be a reduced operad, and let $A$ be a cooperad up to homotopy such
that $\pi_A:T(A[-1]),\pd)\longrightarrow Q$ is a quasi isomorphism,
and let $\phi:Q\longrightarrow P$ be an operad under $Q$. Then
$L_A(P)$ is quasi isomorphic to any complex as in Definition
\ref{Df:Opcohom}.
\end{Tm}
\begin{Pf}
Let $C$ be a cooperad such that $B^*C$ is a cofibrant replacement for
$Q$. We construct an
$L_\infty$-morphism $\hat{\pmb{\oplus}}_SP^A \longrightarrow
\hat{\pmb{\oplus}}_SP^C$ and show that the the induced map in spectral
sequences is an isomorphism in $E_k$ for $k\geq 1$. 

Let $(TA[-1],\pd)$ be an arbitrary cofibrant replacement for $Q$. Then
the solid diagram
\[
\xymatrix{ 0 \ar@{>->}[r]\ar@{>->}[d] & T(C[-1]) \ar@{->>}^{\pi_C}[d] \\
T(A[-1])\ar@{.>}[ur] \ar[r]^{\pi_A} & Q,}
\]
where left arrow is a cofibration and the right arrow is a
acyclic fibration admits a lift (dotted arrow) compatible with
$\pd$. Since bottom and right 
arrow are weak equivalences, so is the lift. This defines a quasi
isomorphism of cooperads up to homotopy from $A$ to $C$ and thus a quasi
isomorphism of $L_\infty$-algebras $L_C(P) \longrightarrow L_A(P)$.

Since the differential on $L_A(P)$ is defined as
the perturbed differential with respect to the solution of the
Maurer-Cartan equation corresponding to the composition 
$(T'(A[-1]),\pd)\longrightarrow Q \longrightarrow P$,
we can again apply Lemma \ref{Lm:perturbmap} and Proposition
\ref{Pp:specseq}.
\end{Pf}

Remember that the \note{cotangent cohomology} $H^{**}(Q;P)$ defined by Markl
\cite{Mar:Model} is computed using a special kind of
minimal model, a so called bigraded minimal model (do not confuse this
notion with the bigraded $L_\infty$-algebras we defined before). Let $Q$ be a
reduced operad, let $\phi:Q\longrightarrow P$ be an operad under 
$Q$, and let $T(M[-1])$ be a bigraded minimal model for $Q$. Then the
cotangent complex can be written as $(\pmb{\oplus}_SP^M,d_\phi)$,
where $d_\phi$ is the internal differential of the
$\phi\circ\pi_M$-perturbed $L_\infty$-algebra $L_A(P)$. 

\begin{Cr}\label{Tm:cotan}
Let $Q$ be a reduced operad, and let $\phi:Q\longrightarrow P$ be an
operad under $Q$. Then
\[
H_Q^p(P) = \bigoplus_q H^{pq}(Q;P).
\]
\end{Cr}
\begin{Pf}
Apply Theorem \ref{Tm:nostrictcohom} to $A=M$, a bigraded minimal model.
\end{Pf}

Let $Q$ be a 1-reduced operad. An operad under $Q$ is an operad $P$ together
with an operad map $\phi:Q\longrightarrow P$. 
This section gives an elegant purely operadic definition of an
$L_\infty$-algebra $L_A(P)$ controlling deformations of operads under
$Q$ (depending on a cofibrant resolution $B^*A$). That is, we have a
deformation theory of operad maps.
Deformations of the map $\phi:Q\longrightarrow P$ correspond to solutions
of the Maurer Cartan equation. There is a natural equivalence reation,
given by a suitable notion of homotopy for coalgebras
(cf. eg. Schlessinger and Stasheff \cite{SchSta:htpy}). The quotient
of the space of solutions by this equivalence relation is invariant
under quasi isomorphism. Therefore the quotient space is
independant of the choice of model for $Q$ in the construction of the
$L_\infty$-algebra. We call this quotient space the \note{space of
deformations} of the map $\phi$. We will not be explicit about the
equivalence relation here, though one might get some intuition from
the next section.

Due to Corollary \ref{Tm:cotan}, this is a more precise statement of Markl's
observation: ``The cotangent cohomology is the best possible
cohomology controlling deformations of operads.''

In Kontsevich's terminology \cite{Kon:Qdef}, the $L_\infty$-algebra
$L_A(P)$, computing the operad cohomology $H^*_Q(P)$, is a formal
pointed manifold with  an odd square-zero vector field controlling
deformations of the operad map $\phi:Q\longrightarrow P$. To phrase
this in terms of formal deformations, pass to $L_A(P)\otimes k[[t]]$,
and consider solutions of Maurer-Cartan of strictly positive degree in
$t$. 

\begin{Ex}\label{Ex:KonSoi}
In the situation described above, let $(T(A[-1]),\pd)$ be a resolution
of $Q$ and let $P=\End_V$. Then the perturbed  $L_\infty$-algebra
$L_A(P)$ of Theorem \ref{Tm:nostrictcohom} defines an
$L_\infty$-algebra. This defines the $L_\infty$-algebra controlling
deformations of $Q$-algebras described by Kontsevich and Soibelman
\cite{KonSoi:Del} by an elegant, purely operadic construction.
\end{Ex}

\subsection{Quadratic Operads}

Let $Q$ be an operad and $\phi:Q\longrightarrow P$ an operad map.
Let $T'(A[-1],\pd)$ be a cofibrant replacement for $Q$.
In Kontsevich's terminology \cite{Kon:Qdef} the $L_\infty$-algebra $L_A(P)$ 
as in Theorem \ref{Tm:nostrictcohom} is a formal pointed manifold with 
an odd square-zero vector field controlling deformations of the 
operad map $\phi:Q\longrightarrow P$. This is illustrated by Theorem
\ref{Tm:OpdMC}, Corollary \ref{Tm:cotan} and Example
\ref{Ex:KonSoi}. This section intends to show that in the case of
quadratic operads our point of view explains a known result of
Balavoine \cite{Bal:Def}.

Let $Q$ be a quadratic operad (cf. Ginzburg and Kapranov
\cite{GinKap:Koszul}) concentrated in degree 0, and let
$\phi:Q\longrightarrow P$ be an operad under $Q$. Let $Q^\bot$ be the dual
cooperad of $Q$ (cf. Getzler and Jones \cite{GetzJon:Opd}). Then we can form
the operad $P^{Q^\bot}$. 
Let $\phi:Q\longrightarrow P$ be an operad under $Q$. We can use the
natural inclusion  $Q^\bot\longrightarrow B^*Q$ to define a
map $\phi^\bot:B^*Q^\bot \longrightarrow P$  of dg operads by 
\begin{equation}\label{eq:compbot}
\phi^\bot:B^*Q^\bot \longrightarrow B^*BQ \longrightarrow Q \longrightarrow
P.
\end{equation}
By slight abuse of notation we denote the $\phi^\bot$-perturbed
$L_\infty$-algebra $\hat{\pmb{\oplus}}_SP^{Q^\bot}$ by $L_{Q^\bot}(P)$.
Balavoine \cite{Bal:Def} shows that the cohomology in
degrees 1 and 2 of $L_{Q^\bot}$ controls formal deformations in the
sense of Gerstenhaber \cite{Gerst:Deform}. The following result shows
that Balavoine's construction gives a nice approximation of
operad cohomology in low degrees, which explains the result.

\begin{Pp}\label{Pp:Baldef}
Let $Q$ be a quadratic operad concentrated in degree 0. Let
$\phi:Q\longrightarrow P$ be an operad under $Q$ and define
$L_{Q^\bot}(P)$ as above. 
The map of complexes $ L_{BQ}(P)\longrightarrow L_{Q^\bot}(P)$ induced
by $ Q^\bot \longrightarrow BQ$ induces an isomorphism 
\[
H^k_Q(P) \longrightarrow H^k(L_{Q^\bot}(P)) \qquad (k\leq 2).
\]
\end{Pp}
\begin{Pf}
By definition of $Q^\bot$, the map $ L_{BQ}(P)\longrightarrow
 L_{Q^\bot}(P)$ (cf. Theorem \ref{Tm:cotan}) is an quasi isomorphism
 in arity $\leq 3$. The result now follows from $Q^\bot(n+1) =
 (Q^\bot)^{-n}$ (cf. Ginzburg and Kapranov \cite{GinKap:Koszul},
 Getzler and Jones \cite{GetzJon:Opd}) together with Lemma
 \ref{Lm:perturbmap} and Proposition \ref{Pp:specseq} restricted to
 degree $\leq 2$. 
\end{Pf}
\begin{Rm}
Thus we have set up a deformation theory of operad maps and shown how
this relates to some known approaches to deformation of algebras over
an operad (Examples \ref{Ex:recover} and \ref{Ex:KonSoi}, and
Proposition \ref{Pp:Baldef}), and cotangent cohomology (Corollary
\ref{Tm:cotan}).
\end{Rm}
\bibliographystyle{plain}
\bibliography{../Hopfalgebras/hopf}
\end{document}